\documentclass{siamltex}
\usepackage{amsfonts}
\usepackage{amssymb}
\usepackage{amsfonts}
\usepackage{mathtools}
\usepackage{color}

\newcommand{\R}{{\mathbb R}}

\newcommand{\bc}{\mathbf{c}}

\newcommand{\bu}{\mathbf{u}}
\newcommand{\bv}{\mathbf{v}}
\newcommand{\bw}{\mathbf{w}}
\newcommand{\be}{\mathbf{e}}

\newcommand{\bg}{\mathbf{g}}

\newcommand{\by}{\mathbf{y}}

\newcommand{\bone}{\mathbf{1}}

\newcommand{\el}{\end{list}}
\newtheorem{algorithm}[theorem]{Algorithm}

\newcommand{\balg}{\begin{algorithm}}
\newcommand{\ealg}{\end{algorithm}}
\newcommand{\bl}{\begin{list}{ \ }{
\leftmargin=.325in}}

\title{Functions and eigenvectors of partially known matrices with applications to 
network analysis}\medskip
\author{
Mohammed Al Mugahwi\thanks{Department of Mathematical Sciences, Kent State University, 
Kent, OH 44242, USA. Email: {\tt malmugah@kent.edu}}
\and
Omar De la Cruz Cabrera\thanks{Department of Mathematical Sciences, Kent State University, 
Kent, OH 44242, USA. Email: {\tt odelacru@kent.edu}}
\and
Silvia Noschese\thanks{Dipartimento di Matematica ``Guido Castelnuovo'', SAPIENZA 
Universit\`a di Roma, P.le A. Moro, 2, I-00185 Roma, Italy. E-mail: 
{\tt noschese@mat.uniroma1.it}}
\and
Lothar Reichel\thanks{Department of Mathematical Sciences, Kent State University, Kent,
OH 44242, USA. E-mail: {\tt reichel@math.kent.edu}}
}
\begin{document}

\maketitle

\begin{abstract}
Matrix functions play an important role in applied mathematics. In network analysis, in
particular, the exponential of the adjacency matrix associated with a network provides
valuable information about connectivity, as well as about the relative importance or 
centrality of nodes. Another popular approach to rank the nodes of a network is to compute 
the left Perron vector of the adjacency matrix for the network. The present article 
addresses the problem of evaluating matrix functions, as well as computing an 
approximation to the left Perron vector, when only some of the columns and/or some of the
rows of the adjacency matrix are known. Applications to network analysis are considered, 
when only some sampled columns and/or rows of the adjacency matrix that defines the 
network are available. A sampling scheme that takes the connectivity of the network into 
account is described.  Computed examples illustrate the performance of the methods 
discussed.
\end{abstract}


%\AMS
%65F10, 65F22
%\endAMS

\section{Introduction}
Many problems in applied mathematics can be formulated and solved with the aid of matrix
functions. This includes the solution of linear discrete ill-posed problems \cite{CR}, 
the solution of time-dependent partial differential equations \cite{DKZ}, and the 
determination of the most important node(s) of a network that is represented 
by a graph and its adjacency matrix \cite{EH,FMRR}. 
Usually, all entries of the adjacency matrix are 
assumed to be known. This paper is concerned with the situation when only some columns, and/or
rows, of the matrix are available. This situation arises, for instance, when one
samples columns, and possibly rows, of a large matrix. We will consider applications in 
network analysis, 
where column and/or row sampling arises naturally in the process of collecting network
data by accessing one node at a time and finding all the other nodes it is connected to.
This is particularly important when
it is too expensive or impractical to collect a full census of all the connections.

A network is represented by a graph $G=\{V,E\}$, which consists of a set
 $V=\{v_j\}_{j=1}^n$ of 
\emph{vertices} or \emph{nodes}, and a set $E=\{e_k\}_{k=1}^m$ of \emph{edges},
the latter being the links between the vertices. Edges may be 
directed, in which case they emerge from a node and end at a node, or undirected. 
Undirected edges are ``two-way streets'' between nodes. For notational convenience and 
ease of discussion, we consider simple (directed or undirected) unweighted graphs $G$ 
without self-loops. Then the 
adjacency matrix $A=[a_{ij}]_{i,j=1}^n\in{\R}^{n\times n}$ associated with the graph $G$ 
has the entry $a_{ij}=1$ if there is a directed edge emerging from vertex $v_i$ and 
ending at vertex $v_j$; if there is an undirected edge between the vertices $v_i$ and 
$v_j$, then $a_{ij}=a_{ji}=1$. Other matrix entries vanish. In particular, the diagonal
entries of $A$ vanish. Typically, $1\le m\ll n^2$, which makes the matrix $A$ sparse. A 
graph is said to be undirected if all its edges are undirected, otherwise the graph is 
directed. The adjacency matrix for an undirected graph is symmetric; for a directed graph 
it is nonsymmetric. Examples of networks include:
\begin{itemize}
\item 
Flight networks, with airports represented by vertices and flights by directed edges.
\item 
Social networking services, such as Facebook and Twitter, with members or accounts 
represented by vertices and interactions between any two accounts by edges.
\end{itemize}
Numerous applications of networks are described in \cite{CEHT,Esbook,Nebook}.

We are concerned with the situation when only some of the nodes and edges of a
graph are known. Each node and its connections to other nodes determine one row
and column of the matrix $A$. Specifically, all edges that point to node $v_i$ 
determine column $i$ of $A$, and all edges that emerge from this node define the 
$i^{\rm th}$ row of $A$. We are interested in
studying properties of networks associated with partially known adjacency matrices.

An important task in network analysis is to determine which vertices of an associated 
graph are the most important ones by measuring how well-connected they are to other 
vertices of the graph. This kind of importance measure often is referred to as a 
\emph{centrality measure}. The choice of a suitable centrality measure depends on what 
the graph is modeling. All commonly used centrality measures ignore intrinsic properties 
of the vertices, and provide information about their importance within the graph just by 
using connectivity information.

A simple approach to measure the centrality of a vertex $v_j$ in a directed graph is to
count the number of edges that point to it. This number is known as the \emph{indegree} of
$v_j$. Similarly, the \emph{outdegree} of $v_j$ is the number of edges that emerge from 
this vertex. For undirected graphs, the \emph{degree} of a vertex is the number of edges
that ``touch'' it.  However, this approach to measure the centrality of a vertex often is
unsatisfactory, because it ignores the importance of the vertices that $v_j$ is connected 
to. Here we consider the computation of certain centrality indices quantifying the 
``importance'' of a vertex on the basis of the importance of its neighbors, according to 
different criteria of propagation of the vertex importance. Such centrality indices are 
based on matrix functions of the adjacency matrix of the graph, and are usually called 
spectral centrality indices. In particular, we focus on the Katz index and the subgraph 
centrality index. Moreover, we also consider eigenvector centrality, that is, the Perron 
eigenvector of the adjacency matrix.

%This shortcoming has prompted the introduction of several centrality measures that are
%based on the evaluation of matrix functions at the adjacency matrix $A$ of $G$; see, e.g.,
%\cite{EH} for a nice introduction. 

To discuss measures determined by matrix functions, we need the notion of a \emph{walk} in
a graph. A walk of length $k$ is a sequence of $k+1$ vertices 
$v_{i_1},v_{i_2},\ldots,v_{i_{k+1}}$ and a sequence of $k$ edges 
$e_{j_1},e_{j_2},\ldots,e_{j_k}$, such that $e_{j_\ell}$ points from $v_{i_\ell}$ to 
$v_{i_{\ell+1}}$ for $\ell=1,2,\ldots,k$. The vertices and edges of a walk do not have to
be distinct. 
It is a well known fact that $[A^k]_{ij}$, i.e., the $(ij)^{\rm th}$
entry of $A^k$, yields the number of walks of length $k$ starting at node $v_i$ and ending
at node $v_j$. Thus, a matrix function evaluated at the adjacency matrix $A$, defined by a
power series $\sum_{k=0}^\infty \alpha_k A^k$ with nonnegative coefficients, can be 
interpreted as containing weighted
sums of walk counts, with weights depending on the length of the walk. Unless $A$ is 
nilpotent (i.e., the graph is directed and contains no cycles), convergence of the power
series requires
that the coefficients $\alpha_k$ converge to zero; this corresponds well with the
intuitively natural requirement that long walks be given less weight than short walks
(which is the case in \eqref{matfun1a} and \eqref{matfun1b} below).

Commonly used matrix functions for measuring the centrality of the vertices 
of a graph are the exponential function $\exp(\gamma_e A)$ and the resolvent 
$(I-\gamma_r A)^{-1}$, where $\gamma_e$ and $\gamma_r$ are positive user-chosen scaling 
parameters; see, e.g., \cite{EH}. These functions can be defined by their power series 
expansions
\begin{eqnarray}\label{matfun1a}
\exp(\gamma_e A)&=&I+\gamma_e A+\frac{1}{2!}(\gamma_e A)^2+\frac{1}{3!}(\gamma_e A)^3+
\ldots~,\\
\label{matfun1b}
(I-\gamma_r A)^{-1}&=&I+\gamma_r A+(\gamma_r A)^2+(\gamma_r A)^3+\ldots~.
\end{eqnarray}
For the resolvent, the parameter $\gamma_r$ has to be chosen small enough so that the 
power series converges, which is the case when $\gamma_r$ is strictly smaller than 
$1/\rho(A)$, where $\rho(A)$ denotes the spectral radius of $A$. 

Matrix functions $f(A)$, such as \eqref{matfun1a} and \eqref{matfun1b}, define several 
commonly used centrality measures: If $f(A)=\exp(A)$, then $[f(A)\bone]_i$ is called the
\emph{total subgraph communicability} of node $v_i$, while the diagonal matrix entry 
$[f(A)]_{ii}$ is the \emph{subgraph centrality} of node $v_i$; see, e.g., \cite{BK,EH}.
Moreover, if $f(A)=(I-\alpha A)^{-1}$, then $[f(A)\bone]_i$ gives the \emph{Katz index}
of node $v_i$; see, e.g., \cite[Chap. 7]{Nebook}. 

It may be beneficial to complement the centrality measures above by the measures 
$[f(A^T)]_{ii}$ and $[f(A^T)\bone]_i$, $i=1,2,\ldots,n$, when the graph $G$ that defines 
$A$ is directed. Here and below the superscript $^T$ denotes transposition; see, e.g., 
\cite{BK,DLCMR,EH,Esbook} for discussions on centrality measures defined by functions of 
the adjacency matrix.

We are interested in computing useful approximations of the largest diagonal entries of 
$f(A)$, or the largest entry of $f(A)\bone$ or $f(A^T)\bone$, when only $1\le k\ll n$ of
the columns and/or rows of $A$ are known. The need to compute such approximations arises when 
the entire graph $G$ is not completely known, but only a small subset of the columns or 
rows of the adjacency matrix $A$ of $G$ are available. This happens, e.g., when not all
nodes and edges of a graph are known, a situation that is common for large,
complex, real-life networks. The
situation we will consider is when the columns and rows of the adjacency matrix are not
explicitly known, but can be sampled. It is then of considerable interest to investigate 
how the sampling should be carried out, as simple random sampling of columns and possibly rows 
of a large adjacency matrix does not give the best results. We will describe a sampling 
method in Section~\ref{sec2}.  A further reason for our interest in computing 
approximations of functions of a 
large matrix $A$, that only use a few of the columns  and/or rows of the matrix, is that 
the evaluation of these approximations typically is much cheaper than the evaluation of 
functions of $A$. 

Another approach to measure centrality is to compute a left or right eigenvector 
associated with the eigenvalue of largest magnitude of $A$. In many situations the
entries of these eigenvectors live in a one-dimensional invariant subspace, have only 
nonvanishing entries, and can be scaled so that all entries are positive. The so-scaled
eigenvectors are commonly referred to as the left and right Perron vectors for the 
adjacency matrix $A$. The left and right Perron vectors are unique up to scaling provided 
that the adjacency matrix is irreducible or, equivalently, if the associated graph is
strongly connected. The centrality of a node is given by the relative size of its 
associated entry of the (left or right) Perron vector for the adjacency matrix. If the 
$j^{\rm th}$ entry of the, say left, Perron vector is the largest, then $v_j$ is the most 
important vertex of the graph. This approach to determine node importance is known as 
\emph{eigenvector centrality} or \emph{Bonacich centrality}; see, e.g., 
\cite{Bo,Esbook,Nebook} for discussions of this method. We will consider the application 
of this method to partially known adjacency matrices.

This paper is organized as follows. Section \ref{sec2} discusses our sampling method for
determining (partial) knowledge of the graph and its associated adjacency matrix. The
evaluation of matrix functions of adjacency matrices that are only partially known is
considered in Section \ref{sec3}, and Section \ref{sec4} describes how an approximation of
the left Perron vector of $A$ can be computed quite inexpensively by using low-rank 
approximations determined by sampling. A few computed examples are presented in Section 
\ref{sec5}, and concluding remarks can be found in Section \ref{sec6}.

\section{Sampling adjacency matrices}\label{sec2}
Let $\sigma_1\geq\sigma_2\geq\dots\geq\sigma_n\geq0$ be the singular values of a large 
matrix $A\in\R^{n\times n}$ and let, for some $1\leq k\ll n$, $\bu_1,\bu_2,\dots,\bu_k$ 
and $\bv_1,\bv_2,\dots,\bv_k$ be left and right singular (unit) vectors associated with
the $k$ largest singular values. Then the truncated singular value decomposition (TSVD) 
\begin{equation}\label{svd}
A^{(k)}=\sum_{j=1}^{k} \sigma_j \bu_j\bv_j^T,
\end{equation}
furnishes a best approximation of $A$ of rank at most $k$ with respect to the spectral and
Frobenius matrix norms; see, e.g., \cite{TB}. However, the computation of the 
approximation \eqref{svd} may be expensive when $n$ is large and $k$ is of moderate size. 
This limits the applicability of the TSVD-approximant \eqref{svd}. Moreover, the 
evaluation of this approximant requires that all entries of $A$ be explicitly known. 

As mentioned above, we are concerned with the situation when $A$ is an adjacency matrix
for a simple (directed or undirected) unweighted graph without self-loops and that, while 
the whole matrix is not known, we can
sample a (relatively small) number of rows and columns.
Then, approximations different from \eqref{svd} have to be used. This
section discusses methods to sample columns and/or rows of $A$. The 
low-rank 
approximations of $A$ determined in this manner are used in Sections \ref{sec3} and 
\ref{sec4} to compute approximations of spectral node centralities.

In the first step, a random non-vanishing column of $A$ is chosen. Let its index
be $j_1$, and denote the chosen column by $\bc_1$. 
If the columns $\bc_1,\dots,\bc_k$ have been chosen, corresponding to the indices
$j_1,\dots,j_k$, at the next step we pick an index $j_{k+1}$ according to a probability
distribution on $\{1,\dots,n\}$ proportional to $\bc_1+\cdots+\bc_k$. 
Thus, at the $(k+1)^{\rm st}$ step, the probability of choosing column $i$ as the next
sampled column is proportional to the number of edges in the network from node $v_i$ to nodes 
$v_{j_1},\dots,v_{j_k}$.  At each step, if a column has
already been picked, or the new column consists entirely of zeros, this choice is discarded
and the procedure is repeated until a new, nonzero column $\bc_{k+1}$ is obtained. We
denote by $J$ the set of indices of the chosen columns; using MATLAB notation, the matrix 
$A_{(:,J)}$ is made up of the chosen columns of $A$. Another way of describing this sampling
method is that we pick the first vertex at random, and then pick subsequent vertices
randomly using a probability distribution proportional to $\bc_1+\cdots+\bc_k$.

We remark that this scheme for selecting columns can just as easily be used in the 
case when the edges have positive weights (that is, the nonzero entries of $A$ may be
positive numbers other than 1). Also, if a row-sampling scheme is needed, rows of the
adjacency matrix $A$
can be selected similarly by applying the above scheme to the columns of the matrix 
$A^T$; in this case we denote by $I$ the set of row indices. The matrix 
$A_{(I,:)}\in\R^{k\times n}$ contains the selected rows of $A$. By alternating column
and row sampling, sets of columns and rows can be determined simultaneously.

The adaptive cross approximation method (ACA) applied to a matrix $A$ also samples rows 
and columns to obtain an approximation of the whole matrix. In ACA, one uses the fact that 
the rows and columns of $A_{(I,:)}$ and $A_{(:,J)}$ have common entries. These 
entries form the matrix $A_{(I,J)}\in\R^{k\times k}$. When the latter matrix is 
nonsingular, the cross approximation of $A$ is given by 
\begin{equation}\label{cross}
M_k=A_{(:,J)}A_{(J,I)}^{-1}A_{(I,:)};
\end{equation}
see \cite{FVB,GTZ,GTZ2,MRVBV} for details.

Let $\sigma_{k+1}\geq 0$ be the $(k+1)^{\rm st}$ singular value of $A$. Then the matrix 
\eqref{svd} satisfies $\|A-A^{(k)}\|_2=\sigma_{k+1}$, where $\|\cdot\|_2$ denotes the 
spectral norm. Goreinov et al. \cite{GTZ} show that there is a matrix $M_k^*$ of rank $k$, 
determined by cross approximation of $A$, such that
\begin{equation}\label{cvbd}
\|A-M_k^*\|_2={\mathcal O}(\sigma_{k+1} \sqrt{kn}).
\end{equation}
Thus, cross approximation can determine a near-best approximation of $A$ of rank $k$ 
without computing the first $k$ singular values and vectors of $A$.

However, the selection of columns and rows of $A$ so that \eqref{cvbd} holds is
computationally difficult. In their analysis, Goreinov et al. \cite{GTZ2} select sets $I$
and $J$ that give the submatrix $A_{(I,J)}$ maximal ``volume'' (modulus of the 
determinant). It is difficult to compute these index sets in a fast manner. Therefore, 
other methods to select the sets $I$ and $J$ have been proposed; see, e.g., 
\cite{FVB,MRVBV}. They are related to incomplete Gaussian elimination with complete 
pivoting. These methods work well when the matrix $A$ is not very sparse. The adjacency 
matrices of concern in the present paper typically are quite sparse, and we found the 
sampling methods described in \cite{FVB,MRVBV} often to give singular matrices 
$A_{(I,J)}$. This makes the use of adaptive cross approximation difficult. We therefore 
will not use the expression \eqref{cross} in subsequent sections.

\section{Functions of low-rank matrix approximations}\label{sec3}
This section discusses the approximation of functions $f$ of a large matrix
$A\in\mathbb{R}^{n\times n}$ that is only partially known. Specifically, we assume 
that only $1\leq\ell\ll n$ columns of $A$ are available, and we would like to determine an
approximation of $f(A)$. We will tacitly assume that the function $f$ and matrix $A$ are 
such that $f(A)$ is well defined; see, e.g., \cite{GVL,Hi} for several definitions of matrix 
functions. For the purpose of this paper, the definition of a matrix function by its power
series expansion suffices; cf. \eqref{matfun1a} and \eqref{matfun1b}. We first will assume
that the matrix $A$ is nonsymmetric. At the end of this section, we will address the 
situation when $A$ is symmetric.

Let $P\in\mathbb{R}^{n\times n}$ be a permutation matrix such that the known columns of 
the matrix $AP$ have indices $1,2,\ldots,\ell$. Thus, the first columns of $AP$ are 
$\bc_1,\ldots,\bc_\ell$. Let $\widetilde{\bc}_j=P^T\bc_j$ for $1\leq j\leq\ell$. We first 
approximate $P^TAP$ by 
\begin{equation}\label{All}
A_\ell=[\widetilde{\bc}_1,\dots,\widetilde{\bc}_\ell,\underbrace{\mathbf{0},\ldots,\mathbf{0}}_{n-\ell}]
\end{equation}
Thus, 
\[
A_\ell=P^TAP\left[\begin{array}{cc} I_\ell & 0 \\ 0 & 0 \end{array}\right]\approx P^TAP,
\]
and then approximate $f(A)=Pf(P^TAP)P^T$ by 
\begin{equation}\label{fapprox}
f(A)\approx Pf(A_\ell)P^T.
\end{equation}
Hence, it suffices to consider the evaluation of $f$ at an $n\times n$ matrix, whose 
$n-\ell$ last columns vanish. We will tacitly assume that $f(A_\ell)$ is well defined.

The computations simplify when $f(0)=0$. We therefore will consider the functions
\begin{equation}\label{modfun}
f(A_\ell)=\exp(\gamma_e A_\ell)-I\mbox{~~~and~~~}f(A_\ell)=(I-\gamma_r A_\ell)^{-1}-I.
\end{equation}
The subtraction of $I$ in the above expressions generally is of no significance for the 
analysis of networks, because one typically is interested in the relative sizes of the 
diagonal entries of $f(A_\ell)$, or of the entries of the vectors $f(A_\ell)\bone$ or
$f(A_\ell^T)\bone$.

The power series representations of the functions in \eqref{modfun},
\[
f(A_\ell)=c_1A_\ell+c_2A_\ell^2+\ldots~,
\]
show that only the first $\ell$ columns of the matrix $f(A_\ell)$ contain nonvanishing
entries. 

Let $\bv_1$ be a random unit vector (not belonging to ${\rm span}\{\bc_1,\ldots,\bc_\ell\}$). 
Application of $\ell$ steps of the Arnoldi process to
$A_\ell$ with initial vector $\bv_1$, generically, yields the Arnoldi decomposition
\begin{equation}\label{arndec}
A_\ell V_{\ell+1}=V_{\ell+1}H_{\ell+1},
\end{equation}
where $H_{\ell+1}\in\mathbb{R}^{(\ell+1) \times (\ell+1)}$ is an upper Hessenberg matrix
and the matrix $V_{\ell+1}\in\mathbb{R}^{n\times(\ell+1)}$ has orthonormal columns. The
computation of the Arnoldi decomposition \eqref{arndec} requires the evaluation of $\ell$
matrix-vector products with $A_\ell$, which is quite inexpensive since $A_\ell$ has at
most $\ell$ nonvanishing columns. We assume that the decomposition \eqref{arndec} exists. 
This is the generic situation. Breakdown of the Arnoldi process, generically, occurs at 
step $\ell+1$; see Saad \cite[Chapter 6]{Sa} for a thorough discussion of the Arnoldi 
decomposition and its computation.

Introduce the spectral factorization 
\begin{equation}\label{Hfact}
H_{\ell+1}=S_{\ell+1}\Lambda_{\ell+1}S_{\ell+1}^{-1},
\end{equation}
which we tacitly assume to exist. This is the generic situation. Thus, the matrix 
$\Lambda_{\ell+1}$ is diagonal; its diagonal entries are the eigenvalues of $H_{\ell+1}$.
We may assume that the eigenvalues are ordered by nonincreasing modulus. Then the last
diagonal entry of $\Lambda_{\ell+1}$ vanishes. It follows that the last column of the 
matrix $S_{\ell+1}$ is an eigenvector that is associated with a vanishing eigenvalue. 
There may be other vanishing diagonal entries of $\Lambda_{\ell+1}$ as well, but this 
will not be exploited. The situation when the factorization \eqref{Hfact} does not exist 
can be handled as described by Pozza et al. \cite{PPS}.

We have
\[
A_\ell V_{\ell+1}S_{\ell+1}=V_{\ell+1}S_{\ell+1}\Lambda_{\ell+1}.
\]
The columns of $V_{\ell+1}S_{\ell+1}$ are eigenvectors of $A_\ell$. The last column of
$V_{\ell+1}S_{\ell+1}$ is an eigenvector that is associated with a vanishing eigenvalue.

Let $\bw_j=V_{\ell+1}S_{\ell+1}\be_j$, $j=1,2,\ldots,\ell$, where $\be_j$ denotes the 
$j^{\rm th}$ column of an identity matrix of appropriate order. Then 
\[
S_n=[\bw_1,\dots,\bw_{\ell}, \be_{\ell+1}, \dots,\be_{n}]\in\mathbb{R}^{n\times n}
\]
is an eigenvector matrix of $A_\ell$, and
\[
A_\ell S_{n}=S_{n}\begin{bmatrix}\Lambda_{\ell}& & & \\& 0 &&\\ && \ddots & \\& & & 0
\end{bmatrix},
\]
where $\Lambda_{\ell}$ is the $\ell\times\ell$ leading principal submatrix of 
$\Lambda_{\ell+1}$. Hence,
\begin{eqnarray}
\nonumber
f(A_\ell)&=&S_{n}f\left(\begin{bmatrix}\Lambda_{\ell}& & & \\& 0 &&\\ && \ddots & \\
 & & & 0\end{bmatrix}\right)S_{n}^{-1} \\
 \nonumber\\
\label{fAll}
 &=& S_{n}\begin{bmatrix}f(\lambda_1) & & & &&\\ 
 & \ddots & & &&\\ & & f(\lambda_{\ell})& & &\\& & &0 &&\\ & &&&\ddots  \\
 & & & &&0\end{bmatrix}S_{n}^{-1},
\end{eqnarray}
where we have used the fact that $f(0)=0$.

To evaluate the expression \eqref{fAll}, it remains to determine the first $\ell$ rows of 
$S_n^{-1}$. This can be done with the aid of the Sherman--Morrison--Woodbury formulas
\cite[p.\,65]{GVL}. Define the matrix $W=[\bw_1,\bw_2,\ldots,\bw_\ell]\in\mathbb{R}^{n\times\ell}$
and let $I_{n,\ell}\in\R^{n\times\ell}$ denote the leading $n\times\ell$ principal submatrix 
of the identity matrix $I\in\R^{n\times n}$. Then the first $\ell$ rows of $S_n^{-1}$ are given by 
$[(I_{\ell,n}W)^{-1},0_{\ell,n-\ell}]$, and we can evaluate
\begin{equation}\label{fAl}
f(A_\ell)=W f(\Lambda_\ell) [(I_{\ell,n}W)^{-1},0_{\ell,n-\ell}],
\end{equation}
where $0_{\ell,n-\ell}\in\R^{\ell\times(n-\ell)}$ denotes a matrix with only zero entries.

Our approximation of $f(A)$ is given by $Pf(A_\ell)P^T$. For a large matrix $A$, the 
computationally most expensive part of evaluating this approximation, when the matrix 
$A_\ell$ is available, is the computation of the Arnoldi decomposition \eqref{arndec}, 
which requires ${\mathcal O}(n\ell^2)$ arithmetic floating point operations. 

We remark that for functions such that 
\begin{equation}\label{ftrans}
f(A)=(f(A^T))^T,
\end{equation}
which includes the functions \eqref{matfun1a} and \eqref{matfun1b}, we may instead sample 
rows of $A$, which are columns of $A^T$, to determine an approximation of $f(A)$ using the
same approach as described above. We remark that equation \eqref{ftrans} holds for all 
matrix functions $f(A)$ that stem from a scalar function $f(t)$ for $t$.

We turn to the situation when the matrix $A\in\R^{n\times n}$ is symmetric, and 
assume that $1\leq\ell\ll n$ of its columns are known. Let the permutation matrix $P$ be
the same as above. Then the first $\ell$ rows and columns of the symmetric matrix 
$A_\ell=P^TAP$ are available. Letting $\bv_1$ be a random unit vector and applying $\ell$
steps of the symmetric Lanczos process to $A_\ell$ with initial vector $\bv_1$ gives,
generically, the Lanczos decomposition
\begin{equation}\label{landec}
A_\ell V_{\ell+1}=V_{\ell+1}T_{\ell+1},
\end{equation}
where $T_{\ell+1}\in\mathbb{R}^{(\ell+1) \times (\ell+1)}$ is a symmetric tridiagonal 
matrix and $V_{\ell+1}\in\mathbb{R}^{n\times(\ell+1)}$ has orthonormal columns. The
computation of the decomposition \eqref{landec} requires the evaluation of $\ell$
matrix-vector products with $A_\ell$. We assume $\ell$ is small enough so that the 
decomposition \eqref{landec} exists. Breakdown depends on the choice of $\bv_1$. 
Typically this assumption is satisfied; otherwise the computations can be modified.
Breakdown of the symmetric Lanczos process, generically, occurs at step $\ell+1$. We 
now can derive a representation of $f(A_\ell)$ of the form \eqref{fAll}, making use 
of the spectral factorization of $T_{\ell+1}$. The derivation in the present situation is
analogous to the derivation of $f(A_\ell)$ in \eqref{fAll}, with the difference that the
eigenvector matrix $S_\ell$ can be chosen to be orthogonal.

\section{The computation of an approximate left Perron vector}\label{sec4} 
Let $A\in\R^{n\times n}$ be the adjacency matrix of a strongly connected graph. Then $A$ 
has a unique left Perron vector $\by=[y_1,y_2,\ldots,y_n]^T\in\R^n$ of unit length with 
all entries positive. As mentioned above, the importance of vertex $v_i$ is proportional 
to $y_i$. When the matrix $A$ is nonsymmetric, the left Perron vector measures the 
centrality of the nodes as receivers; the right Perron vector yield the centrality of the 
nodes as transmitters. 

Assume for the moment that the (unmodified) adjacency matrix $A$ is nonsymmetric. We 
would like to 
determine an approximation of the left Perron vector by using a submatrix determined 
by sampling columns and rows as described in Section \ref{sec2}. Let the set $J$ contain
the $\ell$ indices of the sampled columns of $A$. Thus, the matrix 
$A_{(:,J)}\in\R^{n\times n}$ contains the sampled columns. Similarly, applying the 
same column sampling method to $A^T$ gives a set $I$ of $\ell$ indices; the 
matrix $A_{(I,:)}\in\R^{n\times n}$ contains the sampled rows. We will compute an 
approximation of the left Perron vector of $A$ by applying the power method to the matrix
$M_\ell=A_{(:,J)}A_{(I,:)}$, which approximates $A^2$ (without explicitly forming $M_\ell$). 
We instead also could have applied the power method to $A_{(I,:)}A_{(:,J)}$. Since the
matrix $M_\ell$ is not explicitly stored, the latter choice offers no advantage.

Possible nonunicity of the Perron vector and non-convergence of the power method can be 
remedied by adding a matrix $E\in\R^{n\times n}$
to $M_\ell$, where all entries of $E$ are equal to a small parameter $\varepsilon>0$. The
computations with the power method are carried out without explicitly storing the matrix 
$E$ and forming $M_\ell+E$. The iterations with the power method applied to $M_\ell+E$ are
much cheaper than the iterations with the power method applied to $A$, when $\ell\ll n$. 
Moreover, our method does not require the whole matrix $A$ to be explicitly known. In the
computed examples reported in Section \ref{sec5}, we achieved fairly accurate rankings of
the most important nodes without using the matrix $E$ defined above. Moreover, we found 
that only fairly few rows and columns of $A$ were needed to quite accurately determine the
most important nodes in several ``real'' examples. 

When the adjacency matrix $A$ is symmetric, we propose to compute the Perron vector of the
matrix $M_\ell=A_{(:,J)}A_{(J,:)}$, which can be constructed by sampling the 
columns of $A$, only, to construct $A_{(:,J)}$, since $A_{(J,:)}=A_{(:,J)}^T$. 
Notice that for symmetric matrices the right and left Perron vectors are the same. 
   
\section{Computed examples}\label{sec5}
This section illustrates the performance of the methods discussed when applied to the 
ranking of nodes in several ``real'' large networks. All computations were carried out
in MATLAB with standard IEEE754 machine arithmetic on a Microsoft Windows 10 computer 
with CPU Intel(R) Core(TM) i7-8550U @ 1.80GHz,  4 Cores, 8 Logical Processors 
and 16GB of RAM.

\begin{figure}[htb!]
\begin{center}
\includegraphics[scale=0.9]{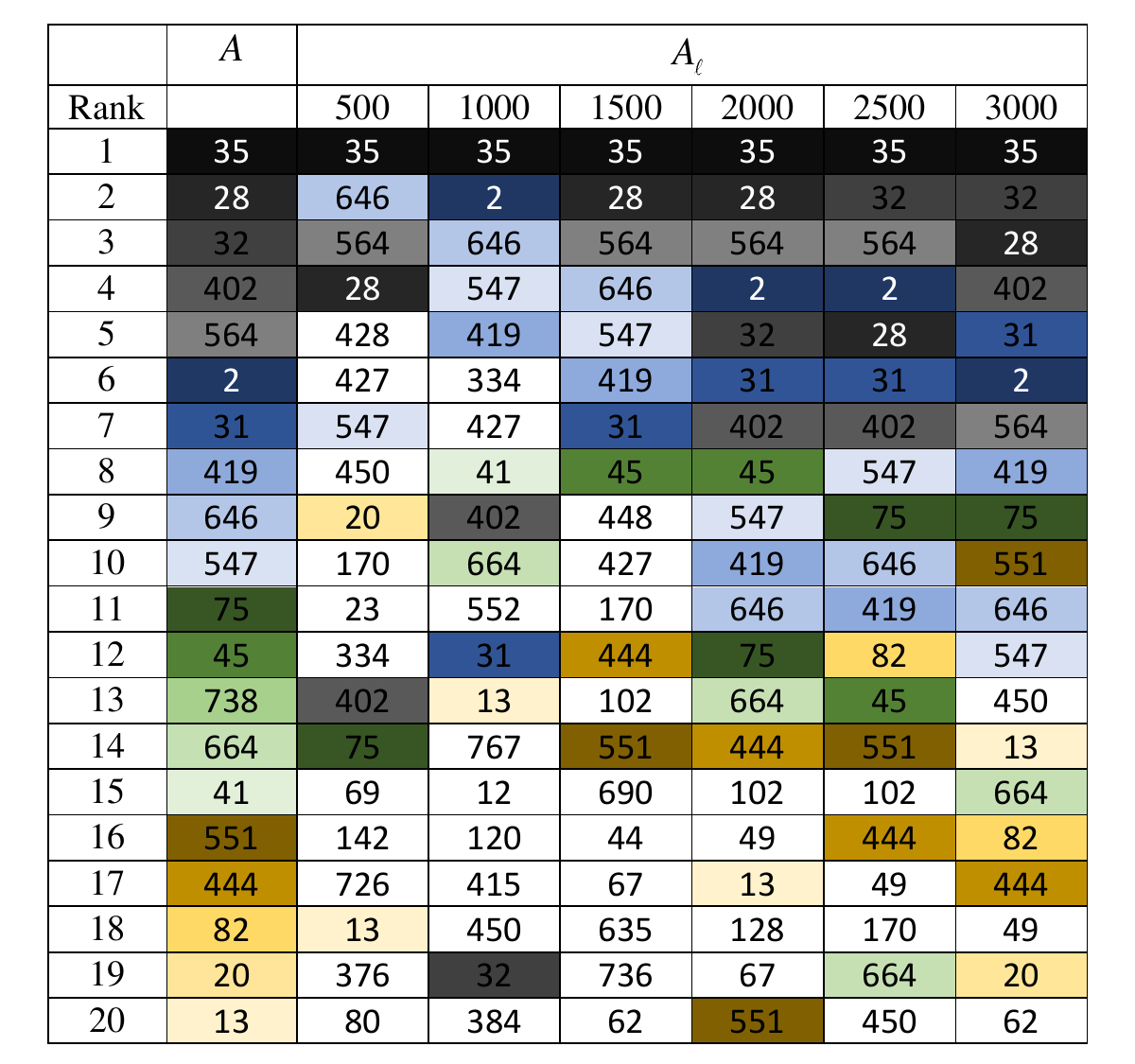}
\end{center}
\caption{soc-Epinions1: The top twenty ranked nodes using the diagonal of $f(A)$ (2nd 
column), and rankings determined by the diagonals of $f(A_\ell)$ for
$\ell\in\{500,1000,1500,2000,2500,3000\}$ for $f(t)=\exp(t)-1$. The columns of $A$ are 
sampled as described in Section~\ref{sec2}.}\label{socFig}
\end{figure}

\begin{figure}[htb!]
\begin{center}
\includegraphics[scale=0.9]{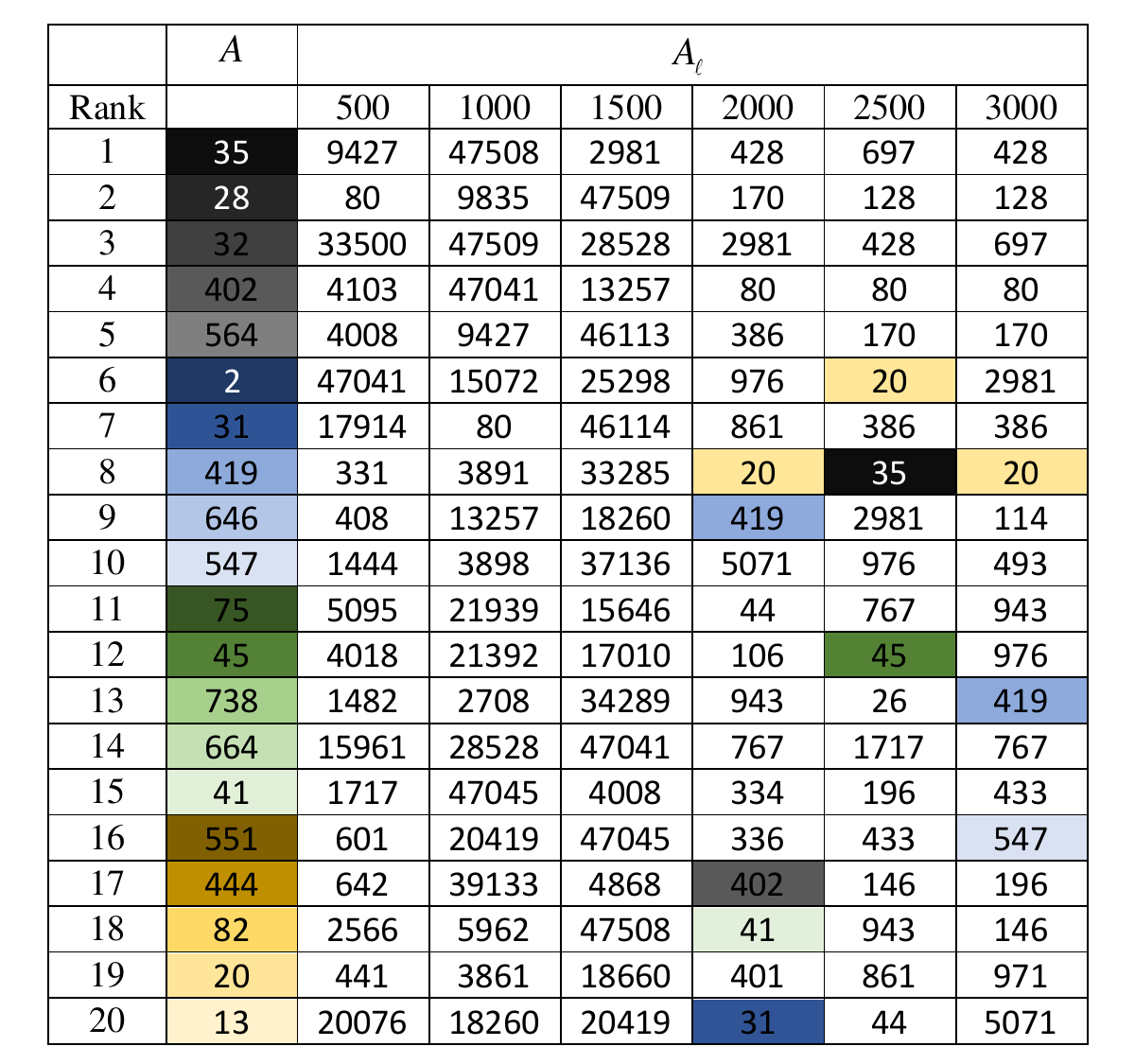}
\end{center}
\caption{soc-Epinions1: The top twenty ranked nodes using the diagonal of $f(A)$ (2nd 
column), and rankings determined by the diagonals of $f(A_\ell)$ for
$\ell\in\{500,1000,1500,2000,2500,3000\}$ for $f(t)=\exp(t)-1$. The columns of $A$ are 
sampled randomly.}\label{socFigrnd}
\end{figure}

\begin{table}[h!]
	\renewcommand{\arraystretch}{1.15}
	\begin{center}
		\begin{tabular}{cccc}
			\hline
			$\ell$ & Mean & Max & Min \\ 
			\hline
			500  & 22.28   & 24.05   & 17.23  \\ 
			1000 & 85.28   & 93.33   & 74.20  \\ 
			1500 & 184.46  & 189.55  & 177.92 \\ 
			2000 & 336.34  & 477.35  & 323.31 \\ 
			2500 & 549.49  & 596.81  & 524.38 \\ 
			3000 & 753.24  & 810.85  & 721.24 \\ 
			\hline \\
		\end{tabular}
	\end{center}
	\caption{soc-Epinions1. Computation time in seconds. Average, max, and min over 
	50 runs.}\label{soc-Epinions1T}
\end{table}

\subsection{soc-Epinions1}\label{sec:soc}
The network of this example is a ``web of trust'' among members of the website 
Epinions.com. This network describes who-trusts-whom. Each user may decide to trust the 
reviews of other users or not. The users are represented by nodes. An edge from node $v_i$
to node $v_j$ indicates that user $i$ trusts user $j$. The network is directed with 75,888 
members (nodes) and 508,837 trust connections (edges) \cite{RAD,SNAP}. We will illustrate 
that one can determine a fairly accurate ranking of the nodes by only using a fairly 
small number of columns of the nonsymmetric adjacency matrix $A\in\R^{n\times n}$ with 
$n=75888$. The node centrality is determined by evaluating 
approximations of the diagonal entries of the matrix function $f(A)=\exp(A)-I$. 

We sample $\ell\ll n$ columns of the adjacency matrix $A$ using the method described in 
Section~\ref{sec2}. The first column, $\bc_1$, is a randomly chosen nonvanishing column of
$A$; the remaining columns are chosen as described in Section~\ref{sec2}. Once the $\ell$
columns of $A$ have been chosen, we evaluate an approximation of $f(A)$ as described in 
Section \ref{sec3}. The rankings obtained are displayed in Figure \ref{socFig}; see below
for a detailed description of this figure. When instead all columns of $A$ are chosen
randomly, then we obtain the rankings shown in Figure \ref{socFigrnd}. Computing times
are reported in Table \ref{soc-Epinions1T}.

The exact ranking of the nodes of the network is difficult to determine due to the large
size of the adjacency matrix. It is problematic to evaluate $f(A)$ both because of the 
large amount of computational arithmetic required, and because of the large storage demand. 
While the matrix $A$ is sparse, and therefore can be stored efficiently using a sparse storage 
format, the matrix $f(A)$ is dense.  In fact, the MATLAB function {\bf expm} cannot be 
applied to evaluate $\exp(A)$ on the computer used for the numerical experiments. Instead,
we apply the Arnoldi process to approximate $f(A)$. Specifically, $k$ steps of the Arnoldi
process applied to $A$ with a random unit initial vector generically gives the Arnoldi 
decomposition 
\begin{equation}\label{arndec2}
AV_k=V_kH_k+\bg_k\be_k^T,
\end{equation}
where the matrix $V_k\in\R^{n\times k}$ has orthonormal columns, $H_k\in\R^{k\times k}$ is
an upper Hessenberg matrix, and the vector $\bg_k\in\R^n$ is orthogonal to the columns of
$V_k$. We then approximate $f(A)$ by $V_kf(H_k)V_k^T$; see, e.g., \cite{BR,DKZ} for 
discussions on the approximation of matrix function using the Arnoldi process. These
computations were carried out for $k=4000$, $k=6000$, $k=8000$, and $k=9000$, and 
rankings ${\rm diag}(V_kf(H_k)V_k^T)$ for these $k$-values were determined. We found
the rankings to converge as $k$ increases. The ranking obtained for $k=9000$ therefore
is considered the ``exact'' ranking. It is shown in the second column of Figure \ref{socFig}.
Subsequent columns of this figure display rankings determined by the diagonal entries of 
$f(A_\ell)$ for $\ell=500$, $1000$, $1500$, $2000$, $2500$, and $3000$, when the columns
of $A$ are sampled by the method of Section \ref{sec2}. Each column shows the top 20 
ranked nodes. To make it easier for a reader to see the rankings, we use $4$ colors, and 
$5$ levels for each color. As we pick $500$ columns of $A$, $9$ of the top $20$ 
ranked nodes are identified, but only the most important node (35) has the correct 
ranking. When $\ell=1000$, the computed ranking improves somewhat. We are able to identify
$11$ out of top $20$ nodes. As we sample more columns of $A$, we obtain improved 
rankings. For $\ell=3000$, we are able to identify $17$ of the $20$ most important nodes, 
and the rankings get closer to the exact ranking. The figure illustrates that useful 
information about node centrality can be determined by sampling many fewer than $n$ 
columns of $A$. Computing times are reported in Table \ref{soc-Epinions1T}.

Figure \ref{socFigrnd} differs from Figure \ref{socFig} in that the columns of the 
matrix $A$ are randomly sampled. Comparing these figures shows the sampling method of 
Section \ref{sec2} to yield rankings that are closer to the ``exact ranking'' of the 
second column for the same number of sampled columns.

\begin{figure}[htb!]
\begin{center}
\includegraphics[scale=0.9]{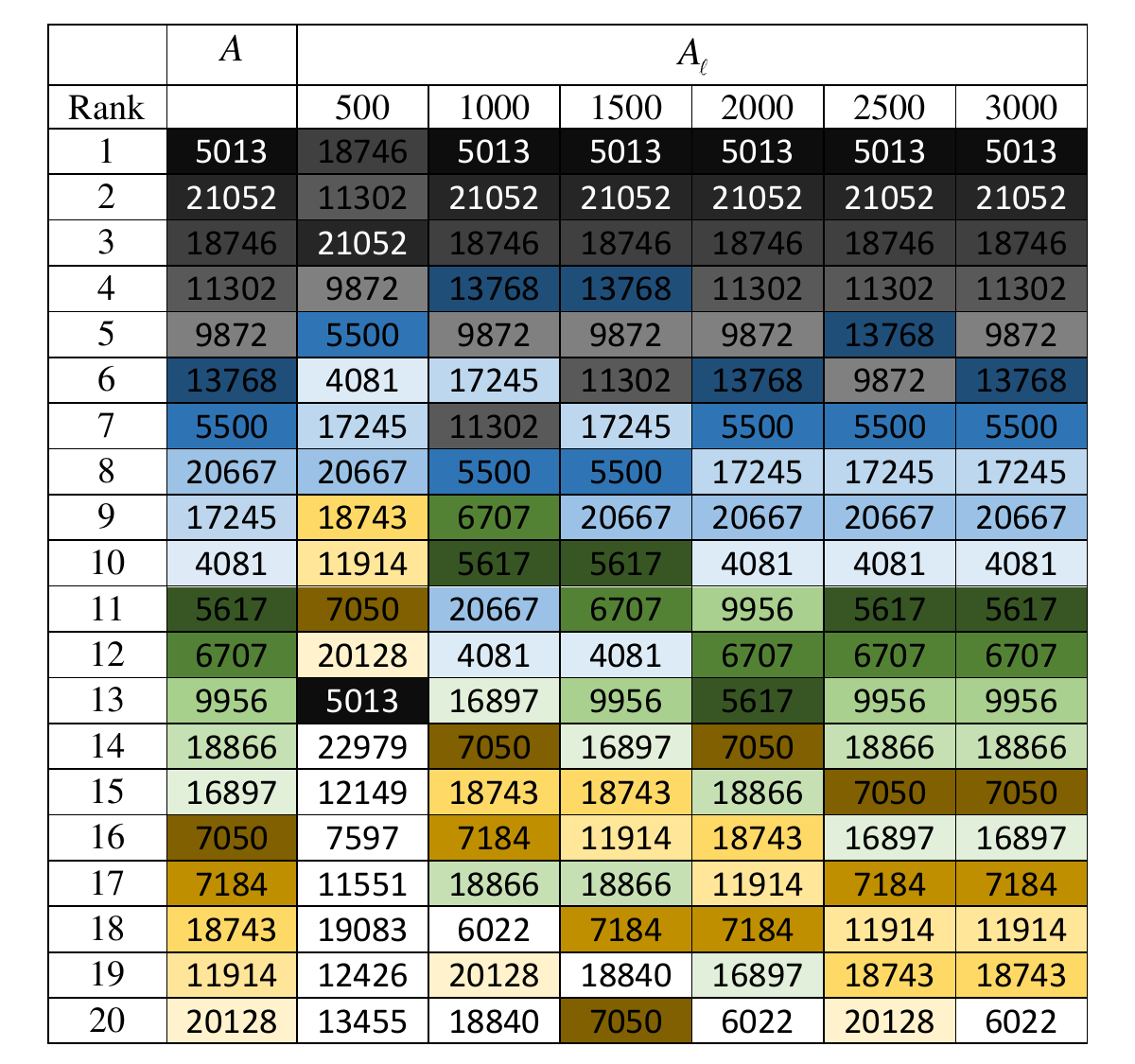}
\end{center}
\caption{ca-CondMat: The top twenty nodes determined by the diagonals of $f(A_\ell)$ for
$\ell\in\{500,1000,1500,2000,2500,3000\}$ for $f(t)=exp(t)-1$. The columns of $A$ are 
sampled as described in Section~\ref{sec2}.}\label{AsFig}
\end{figure}

\begin{figure}[htb!]
\begin{center}
\includegraphics[scale=0.9]{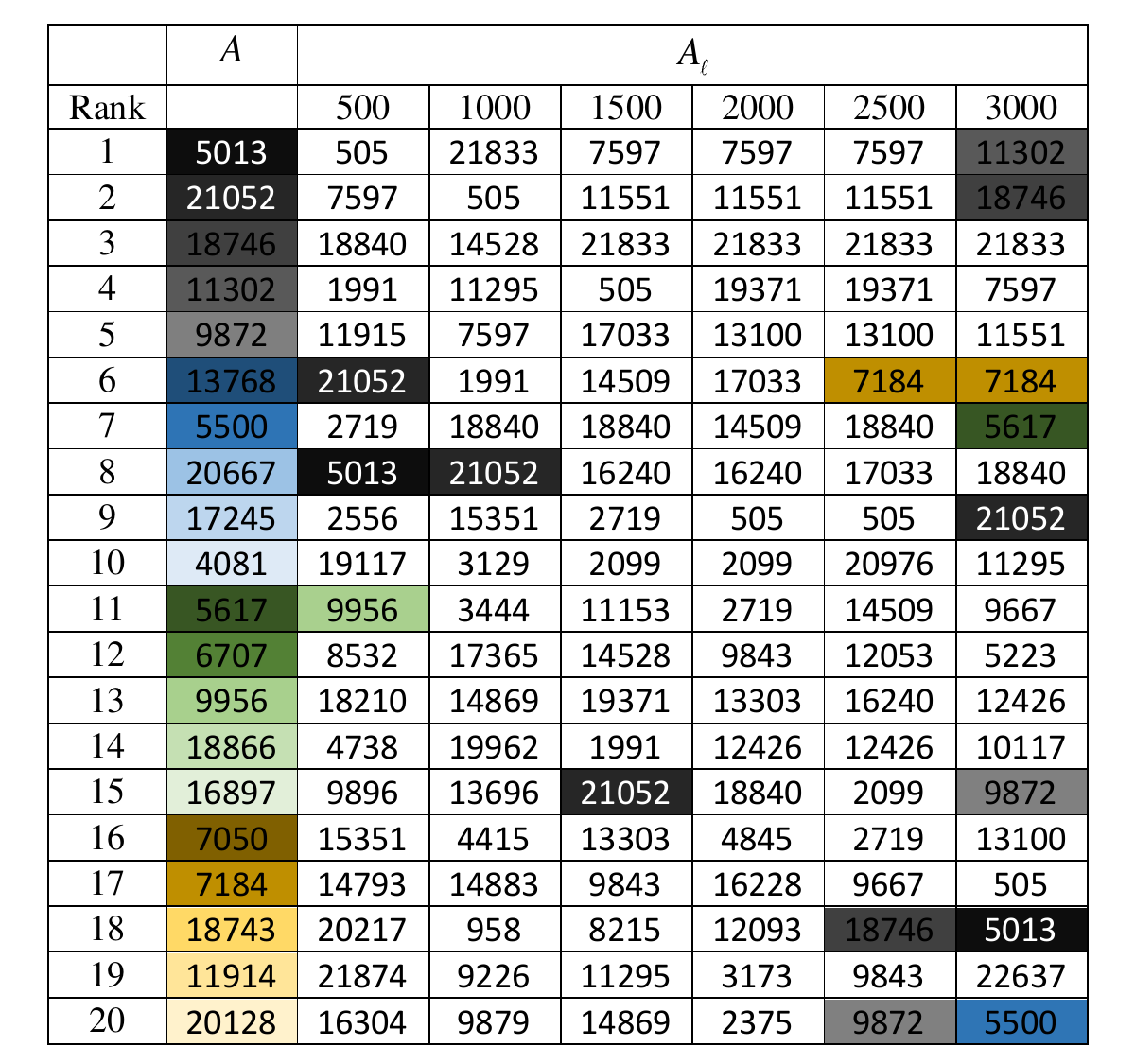}
\end{center}
\caption{ca-CondMat: The top twenty nodes determined by the diagonals of $f(A_\ell)$ for
$\ell\in\{500,1000,1500,2000,2500,3000\}$ for $f(t)=\exp(t)-1$. The columns of $A$ are 
sampled randomly.}\label{AsFigrnd}
\end{figure}

\begin{table}[h!]
	\renewcommand{\arraystretch}{1.15}
	\begin{center}
		\begin{tabular}{cccc}
			\hline
			$\ell$ & Mean & Max & Min \\ 
			\hline
			500  & 0.89  & 1.02  & 0.76  \\ 
			1000 & 2.38  & 3.54  & 2.06  \\ 
			1500 & 4.61  & 7.66  & 3.55  \\ 
			2000 & 7.55  & 12.43 & 5.78  \\ 
			2500 & 10.87 & 19.45 & 8.45  \\ 
			3000 & 15.53 & 37.50 & 11.03 \\ 
			\hline \\
		\end{tabular}
	\end{center}
	\caption{ca-CondMat. Computation time in seconds. Average, max, and min over 
	100 runs. }\label{ca-CondMatT}
\end{table}

\subsection{ca-CondMat}\label{sec:as22} 
This example illustrates the application of the technique of Section \ref{sec3} to a
symmetric partially known matrix. We consider a collaboration network from e-print 
arXiv. The 23,133 nodes of the associated graph represent authors. If author $i$ 
co-authored a paper with author $j$, then the graph has an undirected edge connecting the
nodes $v_i$ and $v_j$. The adjacency matrix $A$ is symmetric with 186,936 non-zero entries 
\cite{LKF,SNAP}. Of the entries, 58 are on the diagonal. Since we are interested in 
graphs without self-loops, we set the latter entries to zero. We use the node centrality 
measure furnished by the diagonal of $f(A)=\exp(A)-I$.  

Figure \ref{AsFig} shows results when using the sampling method described in Section 
\ref{sec2} to choose $\ell$ columns of the adjacency matrix $A$. Due to the symmetry of 
$A$, we also know $\ell$ rows of $A$. The figure compares the ranking of the nodes using 
the diagonal of the matrix $f(A)$ (which is the exact ranking) with the rankings 
determined by the diagonal entries of $f(A_\ell)$ for 
$\ell\in\{500,1000,1500,2000,2500,3000\}$. The figure shows the top 
$20$ ranked nodes determined by each matrix. For $\ell=500$, a couple of the $20$ most 
important nodes can be identified among the first $20$ nodes, but their rankings are incorrect. 
The most important node (5013) is in the $13^{\rm th}$ position, and the second most 
important node (21052) is in the $3^{\rm rd}$ position. Increasing $\ell$ to 1000 yields
more accurate rankings. The  most important nodes, i.e., (5013), (21052), and (18746),
are ranked correctly. Increasing $\ell$ further yields rankings that are closer to the 
``exact'' ranking of the second column. For instance, $\ell=2000$ identifies $19$ of the $20$ most
important nodes, and $8$ of them have the correct rank. The figure  suggests that we 
may gain valuable insight into the ranking of the nodes by using fairly few columns (and
rows) of the adjacency matrix, only.
Computing times are shown in Table \ref{ca-CondMatT}.

Figure \ref{AsFigrnd} differs from Figure \ref{AsFig} in that the columns of the matrix 
$A$ are randomly sampled. Comparing these figures shows that the sampling method of 
Section \ref{sec2} gives rankings that are closer to the ``exact ranking'' of the second 
column for the same number of sampled columns.

\begin{figure}[htb!]
\begin{center}
\includegraphics[scale=1]{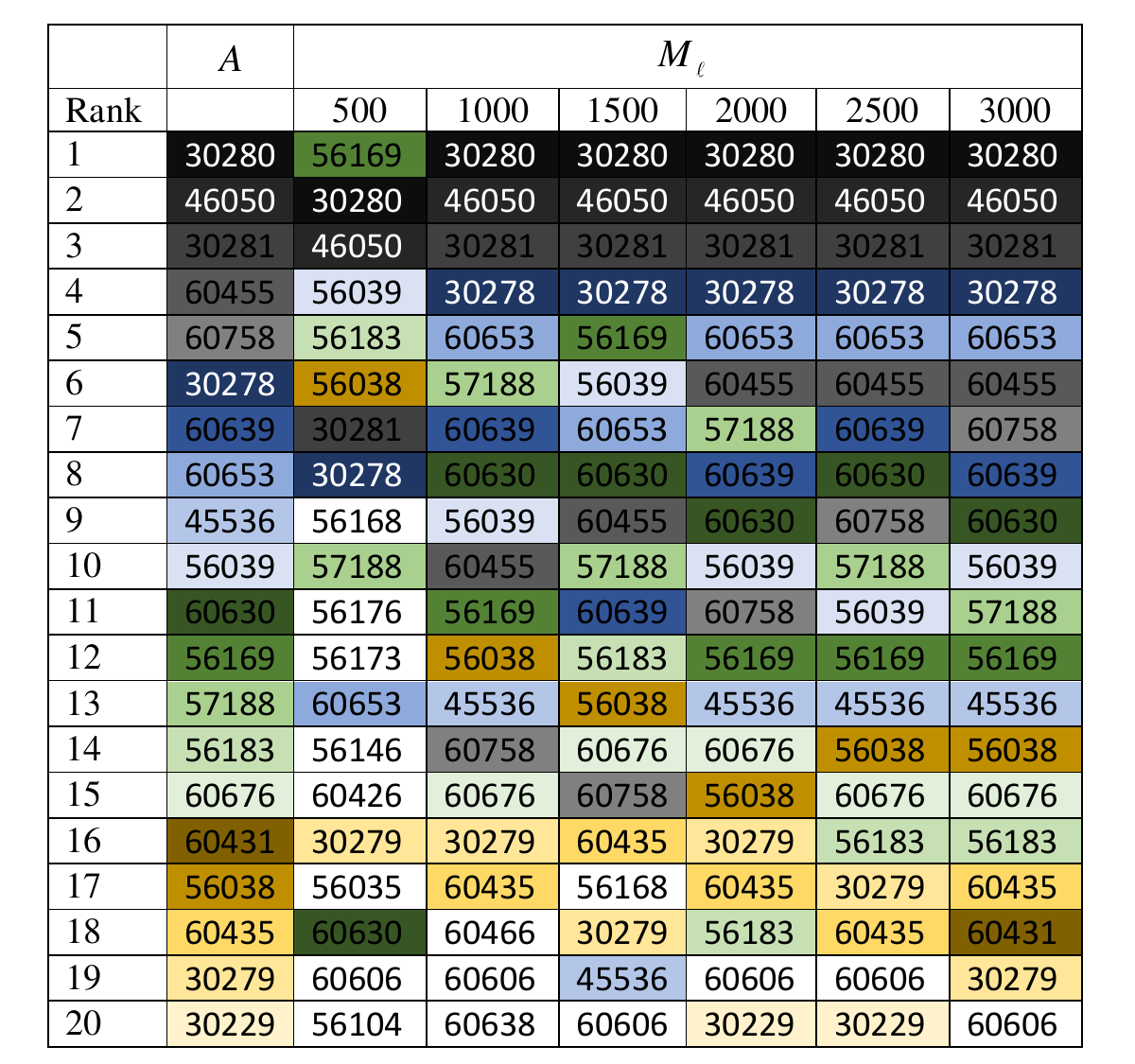}
\end{center}
\caption{Enron: The top $20$ ranked nodes given by the left Perron vector of $A$ and of
$M_\ell=A_{(:,J)}A_{(I,:)}$ for $\ell\in\{500,1000,1500,2000,2500,3000\}$. The columns of 
$A$ are sampled as described in Section~\ref{sec2}.}
\label{EnronTable}
\end{figure}

\begin{figure}[htb!]
\begin{center}
\includegraphics[scale=1]{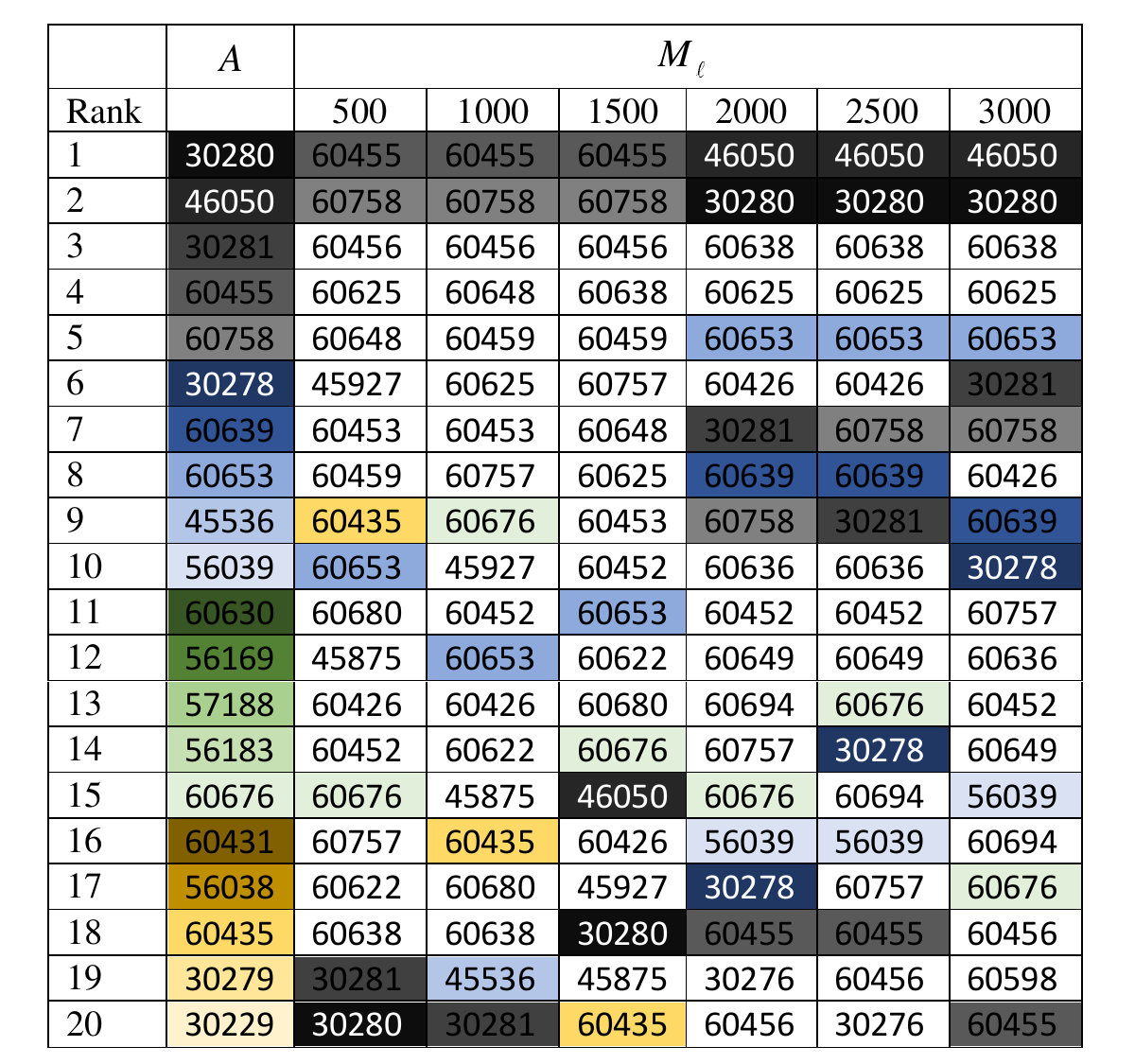}
\end{center}
\caption{Enron: The top $20$ ranked nodes given by the left Perron vector of $A$ and of
$M_\ell=A_{(:,J)}A_{(I,:)}$ for $\ell\in\{500,1000,1500,2000,2500,3000\}$. The columns of
$A$ are sampled randomly.}\label{EnronTablernd}
\end{figure}

\begin{table}[h!]
	\renewcommand{\arraystretch}{1.15}
	\begin{center}
		\begin{tabular}{cccc}
			\hline
			$\ell$ & Mean & Max & Min\\ 
			\hline
			500  & 0.21  & 0.37 & 0.11 \\ 
			1000 & 0.32  & 0.58 & 0.10  \\ 
			1500 & 0.38  & 0.66 & 0.10   \\ 
			2000 & 0.42  & 0.72 & 0.10    \\ 
			2500 & 0.42  & 0.84 & 0.11     \\ 
			3000 & 0.43  & 0.95 & 0.10      \\ 
			\hline \\
		\end{tabular}
	\end{center}
	\caption{Enron.  Computation time in seconds. Average, max, and min over 100 runs.}\label{EnronT}
\end{table}

\subsection{Enron}\label{sec:er}
This example illustrates the application of the method described in Section \ref{sec4} to 
a nonsymmetric adjacency matrix. The network in this example is an e-mail exchange 
network, which represents e-mails (edges) sent between Enron employees (nodes). The
associated graph is unweighted and directed with 69,244 nodes and 276,143 edges, including
1,535 self-loops. We removed the self-loops before running the experiment. This network 
has been studied in \cite{CPAV} and can be found at \cite{SSMC}.

We choose $\ell$ columns of the matrix $A$ as described in Section \ref{sec2} and put
the indices of these columns in the index set $J$. Similarly, we select $\ell$ columns
of the matrix $A^T$. The indices of these rows make up the set $I$. This determines the
matrix $M_\ell=A_{(:,J)}A_{(I,:)}\in\R^{n\times n}$ of rank at most $\ell$. We 
calculate an approximation of a left Perron vector of $A$ by computing a left Perron 
vector of $M_\ell$. The size of the entries of the Perron vectors determines the ranking.

The second column of Figure \ref{EnronTable} shows the ``exact ranking'' determined by a left 
Perron vector of $A$. The remaining columns show the rankings defined by Perron vectors of
$M_\ell$ for $\ell\in\{500,1000,1500,2000,2500,3000\}$ with the sampling of the columns
of $A$ carried out as described in Section~\ref{sec2}. The ranking determined by Perron
vectors of $M_\ell$ gets closer to the exact ranking in the second column as $\ell$ increases.
When $\ell=500$, we are able to identify $12$ out of the $20$ most important nodes, 
but not in the correct order. The three most important nodes have the correct ranking 
for $\ell\geq 1000$. When $\ell\geq 2000$, we almost can identify all the 20 important 
nodes, because node (60606) is actually ranked $21^{\rm st}$. 
Computing times are shown in Table \ref{EnronT}. 

Figure \ref{EnronTablernd} differs from Figure \ref{EnronTable} in that the columns of 
the matrix $A$ are randomly sampled. These figures show that the sampling method of 
Section \ref{sec2} gives rankings that are closer to the ``exact ranking'' of the second 
column for the same number of sampled columns.

\begin{figure}[htb!]
\begin{center}
\includegraphics[scale=1]{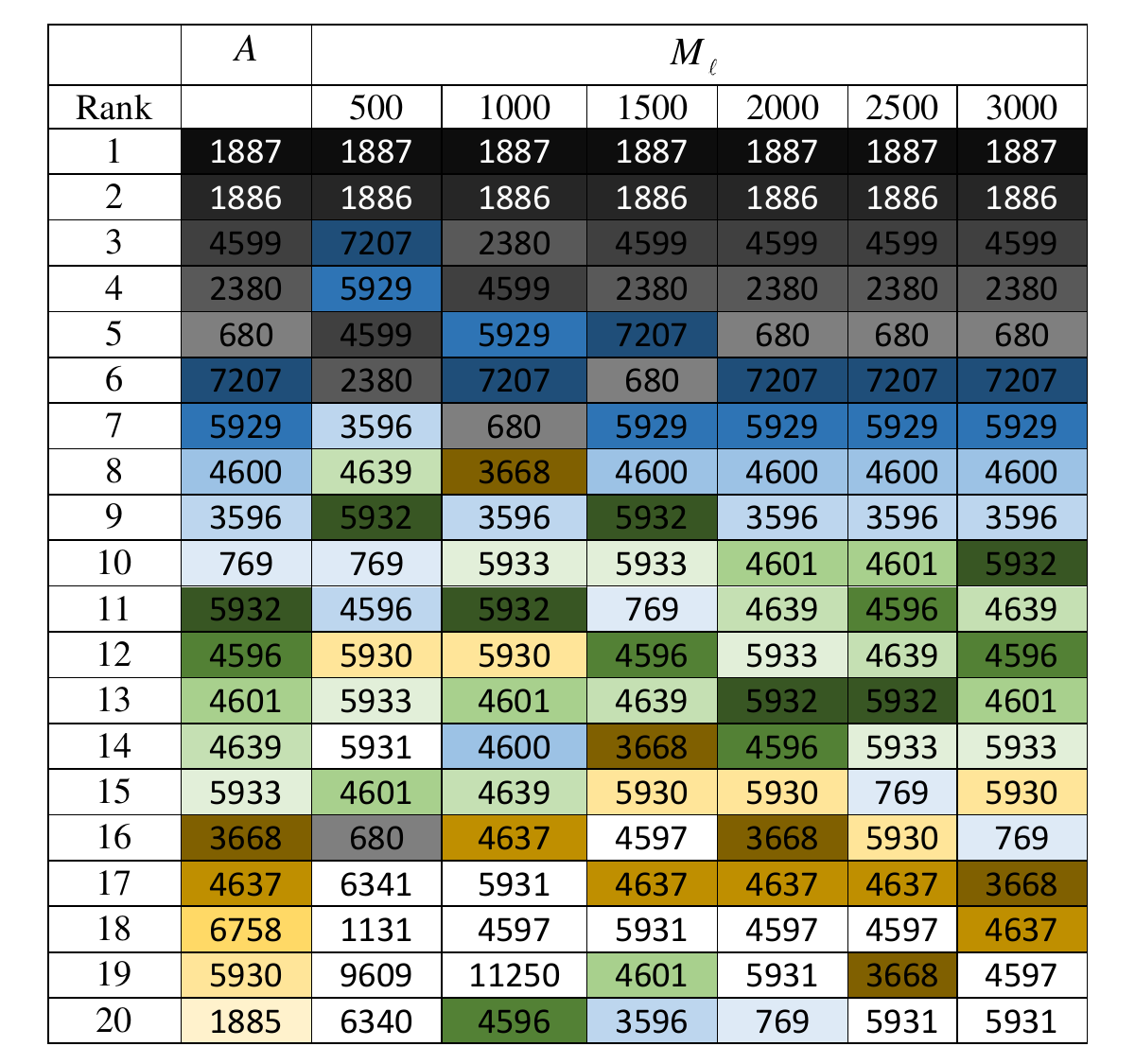}
\end{center}
\caption{Cond-mat-2005: The top $20$ ranked nodes determined by the Perron vectors of $A$
and of $M_\ell=A_{(:,J)}A_{(J,:)}$ for $\ell\in\{500,1000,1500,2000,2500,3000\}$. The 
columns of $A$ are sampled as described in Section~\ref{sec2}.}\label{CondFig}
\end{figure} 

\begin{figure}[htb!]
\begin{center}
\includegraphics[scale=1]{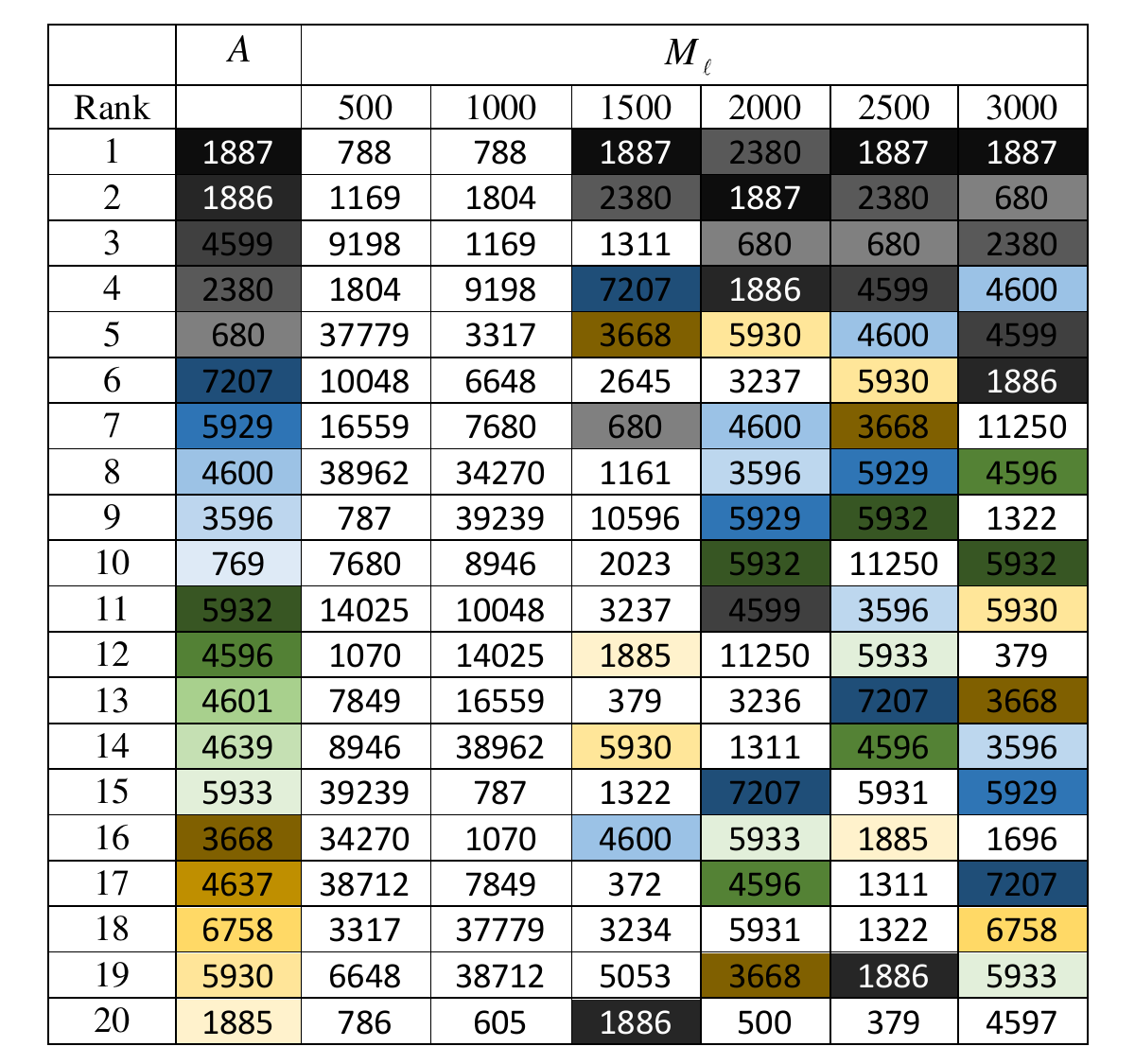}
\end{center}
\caption{Cond-mat-2005: The top $20$ ranked nodes determined by the Perron vectors of $A$
and of $M_\ell=A_{(:,J)}A_{(J,:)}$ for $\ell\in\{500,1000,1500,2000,2500,3000\}$. The 
columns of $A$ are sampled randomly.}\label{CondFigrnd}
\end{figure}

\begin{table}[h!]
	\renewcommand{\arraystretch}{1.15}
	\begin{center}
		\begin{tabular}{cccc}
			\hline
			$\ell$ & Mean & Max & Min   \\ 
			\hline
			500  & 0.08   & 0.13  & 0.03  \\ 
			1000 & 0.13   & 0.17  & 0.02   \\ 
			1500 & 0.16   & 0.24  & 0.03    \\ 
			2000 & 0.19   & 0.26  & 0.03     \\ 
			2500 & 0.21   & 0.25  & 0.03      \\ 
			3000 & 0.23   & 0.28  & 0.03       \\ 
			\hline \\
		\end{tabular}
	\end{center}
	\caption{Cond-mat-2005.  Computation time in seconds. Average, max, and min over 100 runs.}\label{Cond-mat-2005T}
\end{table}

\subsection{Cond-mat-2005}\label{sec:cond}
The network in this example models a
collaboration network of scientists posting preprints in the condensed
matter archive at   www.arxiv.org. It is discussed in \cite{N} and can be found at
\cite{MN}. We use an unweighted version of the network. The associated graph is undirected 
and has 40,421 nodes and 351,382 edges. We use the Perron vector as a centrality measure, 
and compare the node ranking using the Perron vector of $A$ with the ranking determined by 
the 
Perron vector for the matrices $M_\ell=A_{(:,J)}A_{(J,:)}\in\R^{n\times n}$ for several 
$\ell$-values. The matrix $A_{(:,J)}$ is determined as described in Section \ref{sec2},
and $A_{(J,:)}$ is just $A_{(:,J)}^T$. 

Figure \ref{CondFig} shows the (exact) ranking obtained with the Perron vector for $A$ 
(2nd column) and the rankings determined by the Perron vector for $M_\ell$, for 
$\ell\in\{500,1000,1500,2000,2500,3000\}$, when the columns of $A$ are sampled as 
described in Section \ref{sec2}. We compare the ranking of the top $20$ ranked nodes in
these rankings. When $\ell=500$, the two most important nodes are ranked correctly 
by using the Perron vector for $M_{500}$. Moreover, $15$ out of $20$ top ranked nodes are
identified, but their ranking is not correct. For $\ell=2000$, the nine most important
nodes are ranked correctly. 
Computing times are displayed in Table \ref{Cond-mat-2005T}. 

Figure \ref{CondFigrnd} differs from Figure \ref{CondFig} in that the columns of 
the matrix $A$ are randomly sampled. Clearly, the sampling method of Section \ref{sec2} 
gives rankings that are closer to the ``exact ranking'' for the same number of sampled 
columns.

The above examples illustrate that valuable information about the ranking of nodes can 
be gained by sampling columns and rows of the adjacency matrix. The last two examples
determine the left Perron vector. The most popular methods for computing this vector 
for a large adjacency matrix 
are the power method and enhanced variants of the power method that do not require much 
computer storage. These methods, of course, also can be applied to determine the left 
Perron vector of the matrices $M_\ell$. It is outside the scope of the present paper to
compare approaches to efficiently compute the left Perron vector. Extrapolation and other
techniques for accelerating the power method are described in 
\cite{BRZ,BRZ2,BRZ3,CRZT,JS,JS2,WZW}.

In our experience the sampling method described performs well on many ``real'' networks. 
However, one can construct networks for which sampling might not perform well. For instance,
let $G$ be an undirected graph made up of two large clusters with many edges between vertices
in the same cluster, but only one edge between the clusters. The latter edge may be 
difficult to detect by sampling and the sampling method. The method therefore might only
give results for edges in one of the clusters. We are presently investigating how the 
performance of the sampling method can be quantified. We would like to mention that the 
sampling method can be used to study various quantities of interest in network analysis,
such as the total communicability \cite{BK}.

\section{Conclusion}\label{sec6}
In this work we have described novel methods for analyzing large networks in situations
when not all of the adjacency matrix is available. This was done by evaluating matrix
functions or computing approximations of the Perron vector of partially known matrices. 
In the computed examples, we considered the situation when only fairly small subsets of 
columns, or of rows, or both, are known.

There are two distinct advantages to the approaches developed here:
\begin{enumerate}
\item They are computationally much cheaper than the evaluation of matrix functions or 
the computation of the Perron vector of the entire matrix when the adjacency matrix is
large. 
\item The methods described correspond to a compelling sampling strategy when obtaining
the full adjacency information of a network is prohibitively costly. In many realistic
scenarios, the easiest way to collect information about a network is to access nodes (e.g.,
individuals) and interrogating them about the other nodes they are connected to. This
version of sequential sampling is described in Section~\ref{sec2}.
\end{enumerate}

Finally, in order to illustrate the feasibility of our techniques,
we have shown how to approximate well-known node centrality measures for large networks, 
obtaining quite good approximate node rankings, by using only a few columns and rows of 
the underlying adjacency matrix.

\section*{Acknowledgement}
The authors would like to thank Giuseppe Rodriguez and the anonymous referees for 
comments and suggestions.


\begin{thebibliography}{99}
%\bibitem{Be}
%M. Bebendorf, Approximation of boundary element matrices, Numer. Math., 86 (2000), pp. 
%565--589.
\bibitem{BR}
B. Beckermann and L. Reichel, Error estimation and evaluation of matrix functions via the
Faber transform, SIAM J. Numer. Anal., 47 (2009), pp. 3849--3883.
\bibitem{BK}
M. Benzi and C. Klymko, Total communicability as a centrality measure, J. Complex 
Networks, 1 (2013), pp. 124--149.
\bibitem{Bo}
P. F. Bonacich, Power and centrality: A family of measures, Am. J. Sociol., 92 (1987),
pp. 1170--1182.
\bibitem{BRZ}
C. Brezinski and M. Redivo--Zaglia, Rational extrapolation for the PageRank vector,
Math. Comp., 77 (2008), pp. 1585--1598.
\bibitem{BRZ2}
C. Brezinski and M. Redivo--Zaglia, The simplified topological $\varepsilon$-algorithms 
for accelerating sequences in a vector space, SIAM J. Sci. Comput., 36 (2014), pp. 
A2227--A2247.
\bibitem{BRZ3}
C. Brezinski and M. Redivo--Zaglia, The genesis and early developments of Aitken’s 
process, Shanks’ transformation, the $\varepsilon$-algorithm, and related fixed point
methods, Numer. Algorithms, 80 (2019), pp. 11--133.
%\bibitem{BP}
%S. Brin and L. Page, The anatomy of a large-scale hypertextual web search engine,
%Comput. Networks ISDN Systems, 30 (1998), pp. 107--117.
\bibitem{CR}
D. Calvetti and L. Reichel, Lanczos-based exponential filtering for discrete ill-posed 
problems, Numer. Algorithms, 29 (2002), pp. 45--65.
\bibitem{CRZT}
S. Cipolla, M. Redivo-Zaglia, and F. Tudisco, Shifted and extrapolated power methods for
tensor $\ell^p$-eigenpairs, Electron. Trans. Numer. Anal., 53 (2020), pp. 1--27.
\bibitem{CEHT}
J. J. Crofts, E. Estrada, D. J. Higham, and A. Taylor, Mapping directed networks,
Electron. Trans. Numer. Anal., 37 (2010), pp. 337--350.
\bibitem{CPAV}
A. Cruciani, D. Pasquini, G. Amati, and P. Vocca, About Graph Index Compression Techniques, 
Proceedings of the 10th Italian Information Retrieval Workshop (IIR-2019), Padova, Italy, 
September 16-18, 2019, CEUR-WS.org/Vol-2441/paper23.pdf.
\bibitem{DLCMR}
O. De la Cruz Cabrera, M. Matar, and L. Reichel, Analysis of directed networks via the 
matrix exponential, J. Comput. Appl. Math., 355 (2019), pp. 182--192.
\bibitem{DKZ}
V. Druskin, L. Knizhnerman, and M. Zaslavsky, Solution of large scale evolutionary 
problems using rational Krylov subspaces with optimized shifts, SIAM J. Sci. Comput., 31
(2009), pp. 3760--3780.
\bibitem{EH}
E. Estrada and D. J. Higham, Network properties revealed through matrix functions,
SIAM Rev., 52 (2010), pp. 696--714.
\bibitem{Esbook}
E. Estrada, The Structure of Complex Networks, Oxford University Press, Oxford, 2012.
\bibitem{FMRR}
C. Fenu, D. Martin, L. Reichel, and G. Rodriguez, Block Gauss and anti-Gauss quadrature
with application to networks, SIAM J. Matrix Anal. Appl., 34 (2013), pp. 1655--1684.
%\bibitem{FMRR2}
%C. Fenu, D. Martin, L. Reichel, and G. Rodriguez,
%Network analysis via partial spectral factorization and Gauss quadrature, 
%SIAM J. Sci Cumput., 35  (2013), pp. A2046--A2068.
\bibitem{FVB}
K. Frederix and M. Van Barel, Solving a large dense linear system by adaptive cross 
approximation, J. Comput. Appl. Math., 234 (2010), pp. 3181--3195.
\bibitem{GVL}
G. H. Golub and C. F. Van Loan, Matrix Computations, $4$th ed., Johns Hopkins University 
Press, Baltimore, 2013.
\bibitem{GTZ}
S. A. Goreinov, E. E. Tyrtyshnikov, and N. L. Zamarashkin, 
A theory of pseudo-skeleton approximation, Linear Algebra Appl., 261 (1997), pp. 1--21.
\bibitem{GTZ2}
S. A. Goreinov, E. E. Tyrtyshnikov, and N. L. Zamarashkin, 
Pseudo-skeleton approximations by matrices of maximal volume, Math. Notes, 62 (1997), 
pp. 515--519.
\bibitem{Hi}
N. J. Higham, Functions of Matrices: Theory and Computation, SIAM, Philadelphia, 2008.
\bibitem{JS}
K. Jbilou and H. Sadok, LU-implementation of the modified minimal polynomial 
extrapolation method, IMA J. Numer. Anal., 19 (1999), pp. 549--561.
\bibitem{JS2}
K. Jbilou and H. Sadok, Vector extrapolation methods. Applications and numerical 
comparison, J. Comput. Appl. Math., 122 (2000), pp. 149--165.
%\bibitem{LM}
%A. N. Langville and C. D. Meyer, Google’s Pagerank and Beyond, Princeton University 
%Press, Princeton, 2006.
\bibitem{LKF}
J. Leskovec, J. Kleinberg, and C. Faloutsos, Graph evaluation: Densification and shrinking 
diameters, ACM Trans. Knowledge Discovery from Data, 1(1) (2007), Art. 2, pp. 1--41.
%\bibitem{MPI}
%Max Plank Institute for Software Systems, Online Social Networks Research, \hfill\break
%http://socialnetworks.mpi-sws.org/data-wosn2009.html.
\bibitem{MRVBV}
T. Mach, L. Reichel, M. Van Barel, and R. Vandebril, Adaptive cross approximation for 
ill-posed problems, J. Comput. Appl. Math., 303 (2016), pp. 206--217.
\bibitem{MN}
M. E. J. Newman, Network Data, http://www-personal.umich.edu/~mejn/netdata/.
\bibitem{Nebook}
M. E. J. Newman, Networks: An Introduction, Oxford University Press, Oxford, 2010.
\bibitem{N}
M. E. J. Newman, The structure of scientific collaboration networks, 
Proc. Natl. Acad. Sci. USA, 98 (2001), pp. 404--409.
\bibitem{PPS}
S. Pozza, M. S. Prani\'c, and A. Strako\v{s}, The Lanczos algorithm and complex Gauss 
quadrature, Electron. Trans. Numer. Anal., 48 (2018), pp. 362--372.
\bibitem{RAD}
M. Richardson, R. Agrawal, and P. Domingos, Trust management for the semantic web, in
The Semantic Web - ISWC 2003, eds. D. Fensel, K. Sycara, and J. Mylopoulos, Lecture Notes
in Computer Science, vol. 2870, Springer, Berlin, pp. 351--368.
\bibitem{Sa}
Y. Saad, Iterative Methods for Sparse Linear Systems, 2nd ed., SIAM, Philadelphia, 2003.
\bibitem{SNAP}
Stanford Large Network Dataset Collection,
http://snap.stanford.edu/data/index.html
\bibitem{SSMC}
SuiteSparse Matrix Collection, https://sparse.tamu.edu.
\bibitem{TB}
L. N. Trefethen and D. Bau III, Numerical Linear Algebra, SIAM, Philadelphia, 1997.  
\bibitem{WZW}
G. Wu, Y. Zhang, and Y. Wei, Accelerating the Arnoldi-type algorithm for the PageRank 
problem and the ProteinRank problem, J. Sci. Comput., 57 (2013), pp. 74--104.
%\bibitem{Ty}
%E. E. Tyrtyshnikov, Incomplete cross approximation in the mosaic-skeleton method, 
%Computing, 64 (2000), pp. 367--380.
\end{thebibliography}
\end{document}